\documentclass[12pt]{article}
\usepackage{amssymb,amsthm,amsmath,latexsym}
\usepackage{url}
\usepackage{color}
\usepackage{comment}
\usepackage{lscape}
\newtheorem{thm}{Theorem}
\newtheorem{prop}[thm]{Proposition}

\theoremstyle{remark}
\newtheorem{rem}[thm]{Remark}
\newcommand{\FF}{\mathbb{F}}
\newcommand{\ZZ}{\mathbb{Z}}

\newcommand{\cD}{\mathcal{D}}
\newcommand{\cP}{\mathcal{P}}
\newcommand{\cB}{\mathcal{B}}
\DeclareMathOperator{\wt}{wt}
\DeclareMathOperator{\Aut}{Aut}

\begin{document}
\title{
Hadamard matrices of orders $60$ and $64$ with automorphisms
of orders $29$ and $31$}

\author{
Makoto Araya\thanks{Department of Computer Science,
Shizuoka University,
Hamamatsu 432--8011, Japan.
email: {\tt araya@inf.shizuoka.ac.jp}},
Masaaki Harada\thanks{
Research Center for Pure and Applied Mathematics,
Graduate School of Information Sciences,
Tohoku University, Sendai 980--8579, Japan.
email: \texttt{mharada@tohoku.ac.jp}}
and
Vladimir D. Tonchev\thanks{
Department of Mathematical Sciences,
Michigan Technological University,
Houghton, MI 49931, USA.
email: \texttt{tonchev@mtu.edu}.}
}

\maketitle

\begin{abstract}
A classification of Hadamard matrices of order $2p+2$ with an automorphism of order $p$
is given for $p=29$ and $31$.
The ternary self-dual codes spanned by the newly found Hadamard matrices of order $60$
with an automorphism of order $29$ are computed, as well as the binary doubly even 
self-dual codes of length $120$ with generator matrices defined by related Hadamard designs. 
Several new ternary near-extremal self-dual codes,
as well as binary near-extremal doubly even  self-dual codes
with previously unknown weight enumerators are found.
\end{abstract}



\section{Introduction}\label{sec:Intro}

We assume familiarity with basic notions
from error-correcting codes, combinatorial designs and Hadamard matrices~\cite{BJL}, \cite{CRC}, \cite{HP} and \cite{ton88}.
Hadamard matrices  appear in many 
research areas of mathematics and practical applications (see e.g., \cite{Seberry}
and~\cite{SeberryYamada}). 

It is known that for every odd prime $p$ there exists a Hadamard matrix of order $2p+2$
with an automorphism of order $p$, found by Paley~\cite{Paley}, and known in the
combinatorial literature as the Paley-Hadamard matrix of type II.
If $p$ is an odd prime such that $p \equiv -1 \pmod 3$ then the Paley-Hadamard matrix of type II
of order $2p+2$ is a generator matrix of a Pless symmetry code~\cite{Pless},
being a ternary self-dual code of length $2p+2$.
In the context of ternary codes, we consider the elements $0,1,-1$ of $\ZZ$ as
the elements $0,1,2$ of the finite field $\FF_3$ of order $3$.
By the Assmus--Mattson theorem~\cite{AM},
the Pless symmetry codes of length $2p+2$ for $p=5, 11, 17, 23$ and $29$ are ternary extremal self-dual codes
that support 
$5$-designs~\cite{Pless}.
The ternary extended quadratic residue codes of length $2p+2$ are ternary extremal self-dual
codes that support $5$-designs for $p=5, 11, 23$ and $29$, and these codes have generator matrices,
being Hadamard matrices of order $2p+2$ with an automorphism of order $p$~\cite{T}.
In 2013, Nebe and Villar~\cite{NV} gave a series of ternary self-dual codes 
of length $2p + 2$ for a prime $p \equiv 5 \pmod{8}$.
The self-dual codes are also ternary extremal self-dual
codes that support $5$-designs for $p=5$ and $29$~\cite{NV}.
Recently, it was shown in~\cite{AHM} that 
these codes found by Nebe and Villar~\cite{NV} also have generator matrices,
being Hadamard matrices of order $2p+2$ 
for $p \equiv 5 \pmod{24}$.
In addition, it was also shown that there are
at least two inequivalent Hadamard matrices in a third ternary extremal self-dual code of length $60$.
This implies that there are at least four inequivalent Hadamard matrices of order $60$ 
formed by codewords of weight $60$ in
the known three ternary extremal self-dual codes of length $60$.
In addition,  each of these  four Hadamard matrices has an automorphism of order $29$.
This motivates us to classify the Hadamard matrices of order $60$ 
with an automorphism of order $29$.

The Paley-Hadamard matrices of type II for $p=3$ and $p=5$ coincide with the unique, up to equivalence,
Hadamard matrices of orders $8$ and $12$, respectively.  The Hadamard matrices of order $2p+2$
with an automorphism of orders $p=7, 11, 13, 17, 19$ and $23$ have been previously classified up to equivalence.
A short summary on these matrices is given  in Section~\ref{Sec:Summary},  along with relevant references.

In this paper, we present the classification of  Hadamard matrices of order $2p+2$ with an automorphism of order $p$
for the next two odd primes, $p=29$ and $31$.
As an application, we compute 
the ternary self-dual codes and self-dual codes over the finite filed $\FF_5$ of order $5$
spanned by the newly found Hadamard matrices of order $60$
with an automorphism of order $29$,
as well as some binary doubly even self-dual codes of length $120$ with generator matrices defined
by related Hadamard designs. In addition, we study the binary self-dual codes defined 
by Hadamard designs arising from Hadamard matrices of order $64$ with an automorphism of order $31$.

This paper is organized as follows.
In Section~\ref{Sec:2}, we give definitions and some known results
about Hadamard matrices, designs and codes used in this paper.
We also review a method from~\cite{DHM}, \cite{tonchev1} and~\cite{tonchev2}
for classifying  Hadamard matrices of order $2p+2$ with an
automorphism of order $p$, where $p$ is an odd prime.

In Section~\ref{Sec:60}, we give the classification of Hadamard $2$-$(59, 29, 14)$ designs
with an automorphism of order $29$ having one fixed point, and show that there are exactly
$531$ non-isomorphic such designs.
Using these designs, we  classify up to equivalence the Hadamard
matrices of order $60$ with an automorphism of order $29$ and show that
the total number of inequivalent Hadamard matrices with this property is  $266$.
%
%
The ternary code and the code over $\FF_5$ 
spanned by the rows of a Hadamard matrix of order $60$ are self-dual.
A computation of the minimum weights of the ternary codes spanned by the $266$ 
inequivalent Hadamard matrices of order $60$ with an automorphism of order $29$
shows that only four of these matrices span a  ternary extremal self-dual code, and these four 
matrices appear in the three known inequivalent  ternary extremal self-dual codes.
Among the remaining codes, several new ternary near-extremal self-dual codes 
with previously unknown weight enumerators are found.
Self-dual codes over $\FF_5$
spanned by the newly found Hadamard matrices of order $60$
with an automorphism of order $29$ are also computed.
In addition, new binary near-extremal doubly even self-dual codes of length $120$
with previously unknown weight enumerators are constructed from some of the
Hadamard $2$-$(59, 29, 14)$ designs with an automorphism of order $29$.

In Section~\ref{Sec:64}, we give the classification of Hadamard $2$-$(63, 31, 15)$ designs
with an automorphism of order $31$ having one fixed point.
There are exactly $826$ non-isomorphic designs with this property.
Using these designs, we classify the Hadamard
matrices of order $64$ with an automorphism of order $31$, and show that there
are $414$ inequivalent matrices with this property.
The extended code of a binary code of length $63$ spanned by the incidence matrix of a Hadamard 
$2$-$(63, 31, 15)$ design is doubly even.  Among the extended codes of the
$826$ non-isomorphic  $2$-$(63, 31, 15)$  designs with an automorphism of order $31$, there are $794$ 
self-dual codes, and among these codes there are $28$ inequivalent extremal doubly even self-dual codes 
of length $64$ with an automorphism of order $31$.



In Section~\ref{Sec:Summary},
we give a summary of the classification of 
Hadamard $2$-$(2p+1,p,(p-1)/2)$ designs with an automorphism of order 
$p$ having one fixed point
and Hadamard matrices of order $2p+2$ with an automorphism of order $p$,
where $p\le 31$ is an odd prime.

All computer calculations in this paper were done by programs in the language \textsc{C}
and programs in \textsc{Magma}~\cite{Magma}.

\section{Preliminaries}\label{Sec:2}

In this section, we give some definitions and known results
about Hadamard matrices, designs and codes used in this paper.
We also review a method from~\cite{DHM}, \cite{tonchev1} and~\cite{tonchev2}
for classifying  Hadamard matrices of order $2p+2$ with an
automorphism of order $p$, where $p$ is an odd prime.

\subsection{Hadamard matrices, designs and codes}

Throughout this paper, $I_n$ denotes the identity matrix of order $n$,
$A^T$ denotes the transpose of a matrix $A$, and
$J$ denotes the all-one matrix of appropriate size.

A \emph{Hadamard} matrix $H$ of order $n$ is an $n \times n$ matrix
whose entries are from $\{ 1,-1 \}$ such that $H H^T = nI_n$.
It is known that
the order $n$ is necessarily $1,2$, or a multiple of $4$.
Two Hadamard matrices $H$ and $K$ are \emph{equivalent}
if there are $(1,-1,0)$-monomial matrices $P$ and $Q$ with $K = PHQ$.
An \emph{automorphism} of a Hadamard matrix $H$ is
an equivalence of $H$ to itself.
The set of all automorphisms of $H$ forms a group under
composition called the \emph{automorphism group} $\Aut(H)$ of $H$.
For orders up to $32$, all inequivalent Hadamard matrices are known
(see~\cite{KT32} for order $32$).

A \emph{$t$-$(v,k,\lambda)$ design} $\cD$ is a pair of a set $\cP$ of $v$ points
and a collection $\cB$ of $k$-element subsets of
$\cP$ (called blocks)
such that every $t$-element subset of $\cP$ is contained in exactly
$\lambda$ blocks.
Often a $t$-$(v,k,\lambda)$ design is simply called a $t$-design.
Two $t$-$(v,k,\lambda)$ designs are \emph{isomorphic} if
there is a bijection between their point sets that maps the blocks
of the first design into the blocks of the second design.
An \emph{automorphism} of a $t$-$(v,k,\lambda)$ design $\cD$ is any
isomorphism of the design
with itself and the set consisting of all automorphisms of
$\cD$ is called the \emph{automorphism group} $\Aut(\cD)$ of $\cD$.
A $t$-design can be represented by its (block by point)
\emph{incidence matrix} $A=(a_{ij})$, where $a_{ij}=1$ if the $i$-th
block contains the $j$-th point and $a_{ij}=0$ otherwise.
The \emph{complementary} (resp.\ \emph{dual}) 
design of a $t$-design with incidence matrix $A$ is 
the design with incidence matrix $J-A$ (resp.\ $A^T$).
The blocks of a $t$-$(v,k,\lambda)$ design $\cD$ which contain
a given point form a $(t-1)$-$(v-1,k-1,\lambda)$ design
called a \emph{derived} design of $\cD$.
A $2$-design is called \emph{symmetric} if the numbers of points and blocks are the same.
A symmetric $2$-$(4t+3,2t+1,t)$ design is called a \emph{Hadamard} $2$-design.

Let $\FF_p=\{0,1,\ldots,p-1\}$ denote the finite field of order $p$, where $p$ is a
prime.
An $[n,k]$ \emph{code} $C$ over $\FF_p$ is a $k$-dimensional vector subspace
of $\FF_p^n$.
In this paper, we consider codes over $\FF_p$ $(p=2,3,5)$ only.
A code over $\FF_2$ and $\FF_3$ are called
\emph{binary} and \emph{ternary}, respectively.
The parameters $n$ and $k$ are called the \emph{length} and \emph{dimension} of $C$,
respectively.
The \emph{weight} $\wt(x)$ of a vector $x\in\FF_p^n$ is 
the number of non-zero components of $x$.
A vector of $C$ is called a \emph{codeword}.
The minimum non-zero weight of all codewords in $C$ is called
the \emph{minimum weight} of $C$.
The \emph{weight enumerator} of $C$ is defined as $\sum_{c \in C} y^{\wt(c)}$.
Two codes $C$ and $C'$ over $\FF_p$  are \emph{equivalent} if there is a
monomial matrix $P$ over $\FF_p$ with $C' = C \cdot P$,
where $C \cdot P = \{ x P\mid  x \in C\}$.

The {\em dual} code $C^{\perp}$ of a code $C$ of length $n$ is defined as
$
C^{\perp}=
\{x \in \FF_p^n \mid x \cdot y = 0 \text{ for all } y \in C\},
$
where $x \cdot y$ is the standard inner product.
A code $C$ is \emph{self-dual} if $C=C^\perp$.
A binary code $C$ is \emph{doubly even} if $\wt(x) \equiv 0\pmod 4$ for all codewords $x\in C$.
A binary doubly even self-dual code  of length $n$ exists if and only if $n$ is divisible by eight.
A ternary self-dual code  of length $n$ exists if and only if $n$ is divisible by four.
It was shown in~\cite{MS-bound} that
the minimum weight $d$ of a binary doubly even self-dual code
(resp.\ ternary self-dual code) of length $n$
is bounded by $d\leq 4 \lfloor n/24 \rfloor+4$ (resp.\ $d\leq 3 \lfloor n/12 \rfloor+3$).
If $d=4\lfloor n/24 \rfloor+4$ (resp.\ $d=4\lfloor n/24 \rfloor$), 
then the binary doubly even  self-dual code is called \emph{extremal} 
(resp.\ \emph{near-extremal}).
If $d=3\lfloor n/12 \rfloor+3$ (resp.\ $d=3\lfloor n/12 \rfloor$), 
then the ternary self-dual code is called \emph{extremal}
(resp.\ \emph{near-extremal}).

\subsection{Hadamard matrices of order $2p+2$ with an automorphism of order $p$}
\label{Sec:M}

A method for classifying  Hadamard matrices of order $2p+2$ with an
automorphism of order $p$, where $p$ is an odd prime with $p>3$,
is given in~\cite{DHM}, \cite{tonchev1} and~\cite{tonchev2}.
Using this method,
a classification of Hadamard matrices of order $2p+2$ with an
automorphism of order $p$
was completed in~\cite{DHM}, \cite{tonchev1} and~\cite{tonchev2}
for $p=13,17,19,23$.
Here we review the method.

Let $p>3$ be an odd prime.
If a Hadamard matrix $H$ of order $2p+2$ has an
automorphism of order $p$, then $H$ is constructed from a Hadamard
$2$-$(2p+1,p,(p-1)/2)$ design with an automorphism of order $p$ having
one fixed point. 
Note that this follows from a well-known connection between
Hadamard matrices and symmetric designs, together with a bound on the
number of fixed points~\cite[p.~82]{dembowski}.
Let  $\cD$ be a Hadamard
$2$-$(2p+1,p,(p-1)/2)$ design  with an automorphism of order $p$ having one fixed point.
Then $\cD$ has an incidence matrix of the form:
\begin{equation}\label{eq:inc}
A  =
 \left( \begin{array}{ccc}
           &          &  1      \\
    M      &     N    &  \vdots \\
           &          &  1  \\
           &          & 0 \\
    P      &    J-Q   & \vdots \\
           &          & 0\\
1 \cdots 1 & 0 \cdots0& 0 \\
          \end{array} \right),
\end{equation}
where $M,N,P$ and $Q$ are $p \times p$ circulant matrices satisfying
\begin{equation}
MJ=NJ=PJ=QJ=\frac{p-1}{2}J. \label{eq:MJ}
\end{equation}
If the circulant matrices $M,N,P$ and $Q$ satisfy~\eqref{eq:MJ},
the matrix $A$ in~\eqref{eq:inc} is
an incidence matrix of a Hadamard $2$-$(2p+1,p,(p-1)/2)$ design if and
only if the following equalities hold:
\begin{align}
MM^T+NN^T&=\frac{p+1}{2}I_p+\frac{p-3}{2}J, \label{eq:MN}\\
PP^T+QQ^T&=\frac{p+1}{2}I_p+\frac{p-3}{2}J, \label{eq:PQ}\\
MM^T+PP^T&=\frac{p+1}{2}I_p+\frac{p-3}{2}J, \label{eq:MP}\\
NN^T+QQ^T&=\frac{p+1}{2}I_p+\frac{p-3}{2}J, \label{eq:NQ}\\
MP^T&=NQ^T. \label{eq:MPNQ}
\end{align}
Define a $(2p+2) \times (2p+2)$ $(1,-1)$-matrix 
\begin{equation}\label{eq:H1}
H=((-1)^{b_{ij}}), 
\end{equation}
where
\begin{equation}\label{eq:H2}
(b_{ij})=\left( \begin{array}{ccccccc}
1      & 1\cdots1 & 1\cdots1 & 1      \\
1      &          &          & 1      \\
\vdots & M        & N        & \vdots \\
1      &          &          & 1      \\
1      &          &          & 0      \\
\vdots & P        & J-Q      & \vdots \\
1      &          &          & 0      \\
1      & 1\cdots1 & 0\cdots0 & 0
\end{array} \right).
\end{equation}
Then $H$ is a Hadamard matrix of order $2p+2$ with an
automorphism of order $p$.

In the next two sections, we use the above method to extend
the classification of Hadamard matrices of order $2p+2$ with an
automorphism of order  $p=29$ and $p=31$.

\section{Hadamard matrices of order $60$ with an automorphism of order $29$}
\label{Sec:60}

Using the method given in Section~\ref{Sec:M},
in this section, we give the classification of Hadamard $2$-$(59, 29, 14)$ designs
with an automorphism of order $29$ having one fixed point.
Using this classification, we give the classification of Hadamard
matrices of order $60$ with an automorphism of order $29$.
We construct ternary self-dual codes and self-dual codes over $\FF_5$
 from the Hadamard
matrices of order $60$.
We also construct  binary doubly even self-dual codes from 
the Hadamard $2$-$(59, 29, 14)$ designs.

\subsection{Hadamard $2$-$(59, 29, 14)$ designs $D_{59,i}$ and 
Hadamard matrices $H_{60,i}$ of order $60$}

As a first step, we completed the classification of Hadamard $2$-$(59, 29, 14)$ designs
with an automorphism of order $29$  having one fixed point.
In order to classify such $2$-designs, we consider incidence matrices of the 
form~\eqref{eq:inc}.
By a program implemented in the language \textsc{C} using functions from
the GNU Scientific Library (GSL),
our exhaustive computer search found all circulant
matrices $M,N,P$ and $Q$ satisfying the conditions~\eqref{eq:MJ}--\eqref{eq:MPNQ}.
In this calculation, we identify
$14$-elements subsets $S_M,S_N,S_P$ and $S_Q$ of $\{1,2,\ldots,29\}$
with the supports of the first rows of $M,N,P$ and $Q$, respectively.
In this way, we found all Hadamard $2$-$(59, 29, 14)$ designs
with an automorphism of order $29$ which need to be checked further for isomorphism.
To test an  isomorphism of  Hadamard $2$-$(59, 29, 14)$ designs,
we employed the algorithm given in~\cite[p.~15, Theorem~1 (b)]{nautymanual}.
This algorithm considers coloured graphs corresponding to Hadamard $2$-$(59, 29, 14)$ designs.
In this calculation, we used \textsc{nauty} software library~\cite{nauty}
for coloured graphs isomorphism testing.
After isomorphism testing, we completed the classification of Hadamard $2$-$(59, 29, 14)$ designs
with an automorphism of order $29$  having one fixed point.

Using this classification, all Hadamard matrices
of order $60$ with an automorphism of order $29$  are obtained as matrices 
of the form~\eqref{eq:H1} and~\eqref{eq:H2},
which need be checked further for equivalence. 
After equivalence testing, we completed the classification of
Hadamard matrices of order $60$ with an automorphism of order $29$.
This was done by using the \textsc{Magma} function \texttt{IsHadamardEquivalent}.
Then we have the following:

\begin{prop}
There are $531$ non-isomorphic Hadamard $2$-$(59, 29, 14)$ designs
with an automorphism of order $29$ having one fixed point.
There are $266$ inequivalent Hadamard matrices of order $60$ with an automorphism of order $29$.
\end{prop}

The incidence matrices of the above $531$ non-isomorphic Hadamard 
$2$-$(59, 29, 14)$ designs $D_{59,i}$ $(i=1,2,\ldots,531)$ and
the above $266$ inequivalent Hadamard matrices $H_{60,i}$ $(i=1,2,\ldots,266)$ of order $60$
can be obtained electronically from~\cite{AHT}.

The automorphism group orders $|\Aut(D_{59,i})|$ of
$D_{59,i}$  were computed with the  \textsc{Magma} function \texttt{AutomorphismGroup},
 and are listed in Table~\ref{Tab:D59}.
The automorphism group orders $|\Aut(H_{60,i})|$ of $H_{60,i}$  were computed
with the \textsc{Magma} function \texttt{HadamardAutomorphismGroup}, and
 are listed in Table~\ref{Tab:60}.

\begin{table}[thb]
\caption{Orders of the automorphism groups of $D_{59,i}$}
\label{Tab:D59}
\centering
\medskip
{\small
\begin{tabular}{l|l}
\noalign{\hrule height1pt}
\multicolumn{1}{c|}{$|\Aut(D_{59,i})|$ }& \multicolumn{1}{c}{$i$}\\
\hline
$29 \cdot 59$ & $24$\\
$2 \cdot 7 \cdot 29$ & $531$\\
$7\cdot 29$ &  $527, 528, 529, 530$\\
$29$ &  others \\

\noalign{\hrule height1pt}
\end{tabular}
}
\end{table}

\begin{table}[thb]
\caption{Orders of the automorphism groups of $H_{60,i}$}
\label{Tab:60}
\centering
\medskip
{\small
\begin{tabular}{l|l}
\noalign{\hrule height1pt}
\multicolumn{1}{c|}{$|\Aut(H_{60,i})|$ }&  \multicolumn{1}{c}{$i$}\\
\hline
$2^3\cdot 3\cdot 5\cdot 29\cdot 59$ & $21$ \\
$2^5\cdot 3\cdot 5\cdot 7\cdot 29$ & $266$ \\
$2^3\cdot 3\cdot 5\cdot 7\cdot 29$ & $264$ \\
$2^2\cdot 7\cdot 29$ & $265$ \\
$2^2\cdot 29$ & others \\
\noalign{\hrule height1pt}
\end{tabular}
}
\end{table}

\subsection{Ternary self-dual codes $C_3(H_{60,i})$}

Let $H$ be a Hadamard matrix of order $60$.
We denote by $C_3(H)$ the ternary code generated by the rows of $H$.
In the context of ternary codes, we consider the elements $0,1,-1$ of $\ZZ$ as
the elements $0,1,2$ of $\FF_3$, respectively.
It is trivial that if $H$ and $H'$ are equivalent Hadamard matrices of order $60$
then $C_3(H) \cong C_3(H')$.
Since $3$ divides $60$ and $3^2$ does not divide $60$, 
by considering the elementary divisors of $H_{60,i}$, 
it follows that the codes $C_3(H_{60,i})$ $(i=1,2,\ldots,266)$ are self-dual (see~\cite[Section~IV]{LPS}).
The minimum weights $d$ of $C_3(H_{60,i})$ were computed with the \textsc{Magma} function \texttt{MinimumWeight},
and are  listed in Table~\ref{Tab:C60}.
Using the \textsc{Magma} function \texttt{IsIsomorphic}, 
we found the following pairs $(i,j)$
of equivalent codes $C_3(H_{60,i}) \cong C_3(H_{60,j})$:
\begin{equation}\label{eq:F3-60-eq}
\begin{split}
(i,j) =&
( 20, 146),
( 57, 179),
( 75, 103),
( 95, 114),
( 96, 190),
\\&
( 129, 194),
( 132, 140),
( 167, 221),
( 223, 241),
( 264, 265),
\end{split}
\end{equation}
and there is no other pair of equivalent codes among $C_3(H_{60,i})$.

\begin{table}[thb]
\caption{Minimum weights of ternary self-dual codes $C_3(H_{60,i})$}
\label{Tab:C60}
\centering
\medskip
{\small
\begin{tabular}{c|l}
\noalign{\hrule height1pt}
$d$ & \multicolumn{1}{c}{$i$}  \\
\hline
$9$  & $239, 256$ \\
$12$ 
& $2, 6, 7, 11, 12, 16, 20, 25, 26, 30, 31, 32, 33, 34, 36, 37, 40, 44, 45, 53, 54$, \\
& $55, 56, 59, 62, 68, 75, 79, 81, 85, 87, 88, 93, 95, 96, 101, 103, 104, 105$, \\
& $106, 111, 112, 113, 114, 115, 117, 120, 123, 125, 127, 129, 130, 132, 134$, \\
& $135, 137, 140, 143, 144, 146, 148, 152, 154, 155, 160, 164, 165, 168, 172$, \\
& $176, 180, 181, 182, 184, 186, 188, 189, 190, 191, 192, 194, 198, 200, 202$, \\
& $203, 206, 212, 216, 217, 218, 219, 222, 223, 226, 227, 230, 234, 235, 241$ \\
& $244, 245, 250, 252, 261, 262, 263$\\
$15$  & others \\
$18$ & $21, 264, 265, 266$ \\
\noalign{\hrule height1pt}
\end{tabular}
}
\end{table}

It is trivial that an automorphism of order $29$ of $H_{60,i}$ 
induces  an automorphism of order $29$ of $C_3(H_{60,i})$.
There are three inequivalent ternary extremal self-dual codes of length $60$
with an automorphism of order $29$~\cite{BCV}.
The three codes are the extended quadratic residue code $QR_{60}$,
the Pless symmetry code $P_{60}$ and the code $NV_{60}$ found by Nebe and Villar~\cite{NV}.
The codes $QR_{60}$ and $P_{60}$ contain a  type I
Paley-Hadamard matrix $H_{P_I}$ and a type II Paley-Hadamard matrix $H_{P_{II}}$, respectively (see~\cite{T}).
The code $NV_{60}$ contains two inequivalent Hadamard matrices $H_{NV_1}$ 
and $H_{NV_2}$ of order $60$~\cite{AHM}.
%
From Table~\ref{Tab:C60}, we have the following:

\begin{prop}
Let $H$ be a Hadamard matrix $H$ of order $60$
with an automorphism of order $29$ such that $C_3(H)$ 
generates a ternary extremal self-dual code.
Then $H$ is equivalent to one of 
the four inequivalent Hadamard matrices 
$H_{P_I}$, $H_{P_{II}}$, $H_{NV_1}$ and $H_{NV_2}$.
\end{prop}

\begin{table}[thbp]
\caption{Ternary near-extremal self-dual codes of length $36$}
\label{Tab:F3}
\centering
\medskip
{\small
\begin{tabular}{c|l}
\noalign{\hrule height1pt}
$\beta$ & \multicolumn{1}{c}{$i$ } \\ 
\hline
$2552$ & $ 13 $ \\
$2668$ & $ 167, 221 $ \\
$2697$ & $ 52, 138 $ \\
$2726$ & $ 51, 100, 248 $ \\
$2755$ & $ 243, 260 $ \\
$2784$ & $ 4 $ \\
$2813$ & $ 163, 201, 257 $ \\
$2842$ & $ 61, 86, 92, 153, 254 $ \\
$2871$ & $ 28, 39, 83, 108, 229 $ \\
$2900$ & $ 38, 84, 118, 122, 142, 177, 224 $ \\
$2929$ & $ 14, 23, 43, 76, 102, 170, 242, 255 $ \\
$2958$ & $ 5, 42, 58, 73, 82, 91, 174, 211, 215 $ \\
$2987$ & $ 24, 57, 77, 107, 133, 139, 179, 197, 208, 210, 238 $ \\
$3016$ & $ 49, 69, 141, 171, 175, 231 $ \\
$3045$ & $ 17, 78, 119, 157, 185, 251 $ \\
$3074$ & $ 8, 15, 35, 67, 90, 116, 128, 145, 193, 195, 213, 228, 236, 246, 253 $ \\
$3103$ & $ 18, 46, 74, 89, 98, 150, 158, 196, 204, 209, 214, 258 $ \\
$3132$ & $ 29, 48, 64, 66, 80, 94, 131, 147, 156, 187, 207, 237, 240 $ \\
$3161$ & $ 3, 9, 47, 70, 259 $ \\
$3190$ & $ 41, 50, 65, 161, 162, 166, 169, 173, 199, 233 $ \\
$3219$ & $ 27, 60, 97, 136, 178, 220 $ \\
$3248$ & $ 63, 99, 121, 149, 159, 247, 249 $ \\
$3277$ & $ 1, 19, 71, 72, 183, 205 $ \\
$3306$ & $ 124, 126, 225 $ \\
$3335$ & $ 109, 110, 232 $ \\
$3364$ & $ 22 $ \\
$3393$ & $ 151 $ \\
$3422$ & $ 10 $ \\
\noalign{\hrule height1pt}
\end{tabular}
}
\end{table}

Recently, some restrictions on the weight enumerators of 
ternary near-extremal self-dual codes of length divisible by $12$ were proved
in~\cite{AH}, namely, 
the weight enumerator of a ternary near-extremal self-dual code of length $60$ is
of the form:
\[
W_{3,60}=
1+\alpha y^{15}
+(3901080 - 15\alpha) y^{18}
+(241456320 + 105\alpha) y^{21}
+\cdots,
\]
where $\alpha =8\beta$ with $\beta \in \{1,2,\ldots,5148\}$~\cite{AH}.
It is known that there is a ternary near-extremal self-dual code of
length $60$ having weight enumerator $W_{3,60}$
for 
\[
\alpha \in 
\{24\beta \mid \beta \in \Gamma_{60,1}\} \cup \{8\beta \mid \beta \in \Gamma_{60,2}\}, 
\]
where $\Gamma_{60,1}$ and $\Gamma_{60,2}$ are listed in~\cite[Tables~22 and 27]{AH}.
It follows from Table~\ref{Tab:C60} that $154$ of the codes $C_3(H_{60,i})$ are near-extremal.
The values $\alpha=8\beta$ in the weight enumerators $W_{3,60}$
for the $154$ near-extremal self-dual codes $C_3(H_{60,i})$
are listed in Table~\ref{Tab:F3}.
This was calculated by the \textsc{Magma} function \texttt{NumberOfWords}.
From Table~\ref{Tab:F3}, we have the following:

\begin{prop}
\label{p2}
There is a ternary near-extremal self-dual code of length $60$ having weight enumerator 
$W_{3,60}$ for $\alpha \in \{8\beta \mid \beta \in \Gamma_{60,3}\}$, where
\[
\Gamma_{60,3}=
\left\{
\begin{array}{l}
2552, 2668, 2697, 2726, 2755, 2784, 2813, 2842, 2871, \\
2900, 2929, 2987, 3016, 3074, 3103, 3132, 3161, 3190, \\
3219, 3248, 3277, 3306, 3335, 3364, 3422
\end{array}
\right\}.
\] 
\end{prop}

We note that no ternary near-extremal codes with weight enumerators given in Proposition~\ref{p2}
were previously known.

\subsection{Binary doubly even self-dual codes $C_2(D_{59,i})$}

Let $A_{59,i}$ be the incidence matrix of a 
Hadamard $2$-$(59, 29, 14)$ design $D_{59,i}$ and let $C_2(D_{59,i})$ be
the binary $[120,60]$ code  generated by the rows of the following matrix:
\[
\left(\begin{array}{ccccccccc}
{} & {} & {}      & {} & {} & 0     & 1  & \cdots &1  \\
{} & {} & {}      & {} & {} & 1     & {} & {}     &{} \\
{} & {} & I_{60} & {} & {} &\vdots & {} &  B_{59,i}    &{} \\
{} & {} & {}      & {} & {} & 1     & {} &{}      &{} \\
\end{array}\right),
\]
where $B_{59,i}=J-A_{59,i}$ is the incidence matrix of the complementary 
design of  $D_{59,i}$.
In the context of binary codes, we consider the elements $0,1$ of $\ZZ$ as
the elements $0,1$ of $\FF_2$, respectively.
It is trivial that if $D$ and $D'$ are isomorphic Hadamard $2$-$(59, 29, 14)$ designs 
then $C_2(D) \cong C_2(D')$.
The binary codes $C_2(D_{59,i})$ are doubly even self-dual codes~\cite{T83}.
The minimum weights $d$ of $C_2(D_{59,i})$  were computed
with the \textsc{Magma} function \texttt{MinimumWeight}, and 
are  listed in Table~\ref{Tab:B120}.
It is trivial that an automorphism of order $29$ of $D_{59,i}$ 
induces  an automorphism of order $29$ of $C_2(D_{59,i})$.
Note that there is no binary extremal doubly even self-dual code
of length $120$ with an automorphism of order $29$~\cite{CKW}.
Also, it is currently not known whether there is a binary extremal doubly even self-dual code of 
length $120$.

\begin{table}[thb]
\caption{Minimum weights of binary doubly even self-dual codes $C_2(D_{59,i})$}
\label{Tab:B120}
\centering
\medskip
{\small
\begin{tabular}{c|l}
\noalign{\hrule height1pt}
$d$ & \multicolumn{1}{c}{$i$ } \\ 
\hline
$ 8$&$531$ \\
$12$
&$7, 8, 19, 20, 123, 124, 125, 126, 131, 132, 133, 134, 167, 168, 185, 186$, \\
&$197, 198, 217, 218, 221, 222, 251, 252, 281, 282, 285, 286, 337, 338$, \\ 
&$371, 372, 397, 398, 407, 408, 417, 418, 421, 422, 497, 498, 511, 512$ \\
$16$& others \\
$20$ &$21, 24, 209, 210, 529, 530$ \\
\noalign{\hrule height1pt}
\end{tabular}
}
\end{table}

From Table~\ref{Tab:B120},
the codes $C_2(D_{59,i})$ ($i =21, 24, 209, 210, 529, 530$) are binary 
near-extremal doubly even self-dual codes of length $120$.
We verified with the \textsc{Magma} function \texttt{NumberOfWords}
that the numbers of codewords of weight $20$ in these codes are
$71862$,
$71862$,
$98484$,
$98484$,
$104052$ and
$104052$, respectively.
If $D_{59,i}$ and $D_{59,j}$ are Hadamard $2$-$(59, 29, 14)$ designs
which are derived designs of  isomorphic $3$-$(60, 30, 14)$ designs,
then $C_2(D_{59,i})$ and $C_2(D_{59,j})$ are equivalent~\cite[Theorem~2]{T83}.
For $(i,j) \in \{(21, 24), (209, 210), (529, 530)\}$, 
we verified with the \textsc{Magma} function \texttt{IsIsomorphic}
that $D_{59,i}$ and $D_{59,j}$ are derived designs of  
isomorphic $3$-$(60, 30, 14)$ designs.
This implies that $C_2(D_{59,i})$ and $C_2(D_{59,j})$ are equivalent.

The weight enumerator of a binary near-extremal doubly even self-dual code of
length $120$ is of the form:
\begin{align*}
W_{2,120}=&
1
+ \alpha y^{20}
+ (39703755 - 20 \alpha )y^{24}
+ (6101289120 + 190 \alpha )y^{28}
\\&
+ (475644139425 - 1140 \alpha )y^{32}
+ (18824510698240 + 4845 \alpha )y^{36}
\\&
+ (397450513031544 - 15504 \alpha )y^{40}
+ \cdots,
\end{align*}
(see~\cite[Section~4]{H20}).
The existence of $528$ inequivalent binary near-extremal doubly even self-dual codes of
length $120$ is known (see~\cite[Proposition~2]{H20} and \cite[Table~1]{YW}).
These codes have difference weight enumerators $W_{2,120}$, 
where $\alpha$ are given in~\cite[Table~3]{H20} and~\cite[Table~1]{YW}.
\begin{rem}
The number of inequivalent binary near-extremal doubly even self-dual codes of
length $120$ constructed in~\cite{H20} was incorrectly  reported as $500$.
We point out that  the correct number is $502$.  
Hence,  $528$ inequivalent binary near-extremal doubly even self-dual codes of
length $120$ were previously known.
\end{rem}

\begin{prop}\label{prop:F2}
\begin{itemize}
\item[\rm (i)]
There is a binary near-extremal doubly even self-dual code of
length $120$ with weight enumerator $W_{2,120}$ for
$\alpha=98484$ and $104052$.
\item[\rm (ii)]
There are at least $530$ inequivalent binary near-extremal doubly even self-dual codes of
length $120$.
\end{itemize}
\end{prop}

We note that no binary near-extremal codes with weight enumerators given in 
Proposition~\ref{prop:F2}~(i)
were previously known.

\subsection{Self-dual codes $C_5(H_{60,i})$ over $\FF_5$}

Let $H$ be a Hadamard matrix of order $60$.
We denote by $C_5(H)$ the code over $\FF_5$ generated by the rows of $H$.
In the context of codes over $\FF_5$, we consider the elements $0,1,-1$ of $\ZZ$ as
the elements $0,1,4$ of $\FF_5$, respectively.
It is trivial that if $H$ and $H'$ are equivalent Hadamard matrices of order $60$
then $C_5(H) \cong C_5(H')$.
In addition, $C_5(H)$ is self-dual~\cite{LPS2}.
The minimum weights $d$ of $C_5(H_{60,i})$ $(i=1,2,\ldots,266)$
were computed with the \textsc{Magma} function \texttt{MinimumWeight},
and are  listed in Table~\ref{Tab:C60F5}.
Let $d(C_5(H_{60,i}))$ and 
$N(C_5(H_{60,i}))$ denote the minimum weight and 
the number of codewords of minimum weight in $C_5(H_{60,i})$, respectively.
For the pair $(i,j)$ such that $d(C_5(H_{60,i}))=d(C_5(H_{60,j}))$
and $N(C_5(H_{60,i}))=N(C_5(H_{60,j}))$,
using the \textsc{Magma} function \texttt{IsIsomorphic}, 
we determined whether $C_5(H_{60,i}) \cong C_5(H_{60,j})$ or not,
where the numbers $N(C_5(H_{60,i}))$ were computed 
with the \textsc{Magma} function \texttt{NumberOfWords}.
Then we found equivalent codes $C_5(H_{60,238}) \cong C_5(H_{60,257})$
and there is no other pair of equivalent codes among $C_5(H_{60,i})$.

Note that the extended quadratic residue code $QR_{60}$ of length $60$ is a self-dual code having
minimum weight $18$.
The largest minimum weight among self-dual codes of length $60$ is
between $18$ and $24$~\cite[Table~9]{GG}.
From Table~\ref{Tab:C60F5},
$C_5(H_{60,21})$ and $C_5(H_{60,266})$ have the largest minimum weight
among currently known self-dual codes of length $60$.
We verified with the \textsc{Magma} function \texttt{NumberOfWords}
that the numbers of codewords of minimum weight in $QR_{60}$, 
$C_5(H_{60,21})$ and $C_5(H_{60,266})$ are
$410640$, $410640$ and $288840$, respectively.
Therefore, there are at least two inequivalent self-dual codes over $\FF_5$ of length $60$
and minimum weight $18$.


\begin{table}[thb]
\caption{Minimum weights of self-dual codes $C_5(H_{60,i})$ over $\FF_5$}
\label{Tab:C60F5}
\centering
\medskip
{\small
\begin{tabular}{c|l}
\noalign{\hrule height1pt}
$d$ & \multicolumn{1}{c}{$i$}  \\
\hline
$12$  & $66, 103, 129, 142, 170, 198, 209, 248, 259$ \\
$14$ & 
$1, 3, 4, 11, 19, 25, 28, 33, 35, 36, 38, 47, 49, 51, 52, 54, 58, 59, 64, $\\
&$68, 72, 76, 78, 79, 81, 84, 85, 87, 89, 90, 91, 94, 96, 97, 105, 108, $\\
&$109, 111, 113, 114, 115, 116, 118, 119, 120, 122, 123, 128, 132,  $\\
&$133, 135, 139, 149, 150, 152,155, 157, 161, 162, 164, 165, 169,  $\\
&$172, 176, 178, 181, 183, 185, 188, 189, 194,195, 196, 199, 200,   $\\
&$202, 206, 208, 211, 212, 218, 221, 222, 223, 224, 225, 226, 228,   $\\
&$233, 235, 236, 238, 240, 244, 245, 246, 252, 253, 256, 257, 260$\\
$15$  & $7, 45, 48, 50, 53, 60, 61, 62, 65, 70, 92, 95, 101, 102, 107, 121, $\\
&$127, 136, 144, 146, 159, 160, 171, 174, 177, 186, 192, 193, 210, $\\
&$215, 219, 220, 230, 239, 243, 247, 255, 263$ \\
$16$ & others  \\
$18$ & $21, 266$ \\
\noalign{\hrule height1pt}
\end{tabular}
}
\end{table}

\section{Hadamard matrices of order $64$ with an automorphism of order $31$}
\label{Sec:64}

Using the method described in Section~\ref{Sec:M},
 we give in this section the classification of Hadamard $2$-$(63, 31, 15)$ designs
with an automorphism of order $31$ having one fixed point.
Using this classification, we give a complete classification of Hadamard
matrices of order $64$ with an automorphism of order $31$.
We construct  binary doubly even codes from 
the Hadamard $2$-$(63, 31, 15)$ designs.

%
%

\subsection{Hadamard $2$-$(63, 31, 15)$ designs $D_{63,i}$ and 
Hadamard matrices $H_{64,i}$ of order $64$}

The approach used in the classification is similar to that
given in the previous section, so
in this section, only results are given.

\begin{prop}
\label{p5}
There are $826$ non-isomorphic Hadamard $2$-$(63, 31, 15)$ designs
with an automorphism of order $31$ having one fixed point.
There are $414$ inequivalent Hadamard matrices of order $64$ with an automorphism of order $31$.
\end{prop}

\begin{rem}
If $p$ is an odd prime dividing the order of the
automorphism group of a Hadamard matrix of order $n\ge 4$,
then either $p$ divides
$n$ or $n-1$, or $p \le n/2-1$~\cite{tonchev1}.
Hence, the largest prime which can divide
the order of the automorphism group of a Hadamard matrix of order $64$ is $31$.
\end{rem}

The incidence matrices of the above 
$826$ non-isomorphic Hadamard $2$-$(63, 31, 15)$ designs $D_{63,i}$
$(i=1,2,\ldots,826)$ and 
the above  $414$ inequivalent Hadamard matrices $H_{64,i}$ $(i=1,2,\ldots,414)$ from Proposition~\ref{p5}
can be obtained electronically from~\cite{AHT}.

The automorphism group orders $|\Aut(D_{63,i})|$ of $D_{63,i}$ and
the automorphism group orders $|\Aut(H_{64,i})|$ of $H_{64,i}$
are listed in Tables~\ref{Tab:D63} and~\ref{Tab:64}, respectively.

\begin{table}[thb]
\caption{Orders of the automorphism groups of $D_{63,i}$}
\label{Tab:D63}
\centering
\medskip
{\small
\begin{tabular}{l|l}
\noalign{\hrule height1pt}
\multicolumn{1}{c|}{$|\Aut(D_{63,i})|$ }&\multicolumn{1}{c}{$i$ }\\
\hline
$2^{15}\cdot 3^4 \cdot5 \cdot7^2 \cdot31  $&$805$\\
$2^5\cdot 5\cdot 31           $
&$790, 791, 793, 797, 806, 807, 808, 809, 810, 812$, \\
&$813, 816, 820, 823$\\
$3 \cdot 5 \cdot31            $&$789, 792, 799$\\
$5  \cdot   31           $&$794, 795, 796, 798, 800, 801, 802, 803, 804, 811$, \\
&$814, 815, 817, 818,
819, 821, 822, 824, 825, 826$\\
$3  \cdot   31           $&$ 765, 766, 767, 768, 769, 770, 771, 772, 779, 780$, \\
& $781, 782, 783, 784, 785, 786$\\
$31                 $&others \\
\noalign{\hrule height1pt}
\end{tabular}
}
\end{table}

\begin{table}[thb]
\caption{Orders of the automorphism groups of $H_{64,i}$}
\label{Tab:64}
\centering
\medskip
{\small
\begin{tabular}{l|l}
\noalign{\hrule height1pt}
\multicolumn{1}{c|}{$|\Aut(H_{64,i})|$ }&\multicolumn{1}{c}{$i$ }\\
  \hline
$2^{28}\cdot 3^4\cdot 5\cdot 7^2\cdot 31$ & $406$\\
$2^{13}\cdot 5\cdot 31$ & $407, 410$\\
$2^{12}\cdot 5\cdot 31$ & $408, 409, 411, 413$\\
$2^8\cdot 3\cdot 5\cdot 31$ & $395, 398, 405$\\
$2^7\cdot 5\cdot 31$ & $396, 397, 399, 403$\\
$2^2\cdot 5\cdot 31$ & $400, 404, 412, 414$\\
$2^2\cdot 3\cdot 31$ & $387, 388, 389, 390, 391, 392, 393, 394$\\
$2\cdot 5\cdot 31$ & $401,402$\\
$2^2\cdot 31$ & others\\
 \noalign{\hrule height1pt}
\end{tabular}
}
\end{table}

\subsection{Binary doubly even codes $C'_2(D_{63,i})$}

Let $A_{63,i}$ be the  incidence matrix of a 
Hadamard $2$-$(63, 31, 15)$  design $D_{63,i}$ and let $C'_2(D_{63,i})$ be
the binary code  generated by the rows of the following matrix:
\[
\left( \begin{array}{ccccccc}
 &     & & 1      \\
 & A_{63,i}   & & \vdots \\
 &     & & 1      \\
\end{array} \right).
\]
It is trivial that if $D$ and $D'$ are isomorphic Hadamard $2$-$(63, 31, 15)$  designs 
then $C'_2(D) \cong C'_2(D')$.
In addition, it is trivial that $C'_2(D_{63,i})$ are doubly even. 
The dimensions  $\dim(C'_2(D_{63,i}))$ were computed
with the \textsc{Magma} function \texttt{Dimension}, and are listed in Table~\ref{Tab:F2-1}.
The minimum weights $d$ of the $794$ doubly even self-dual codes $C'_2(D_{63,i})$ 
were computed with the \textsc{Magma} function
\texttt{MinimumWeight}, and are listed in Table~\ref{Tab:F2-2}.
In addition,  we verified that the $520$ extremal  doubly even self-dual codes $C'_2(D_{63,i})$
are divided into $28$ equivalence classes.
More precisely, each of the $520$ codes is
equivalent to one of the $28$ inequivalent extremal  doubly even self-dual codes $C'_2(D_{63,i})$, 
where
\[
i \in 
\left\{ \begin{array}{l}
1, 2, 3, 4, 6, 8, 9, 10, 11, 14, 15, 17, 18, 20, 24,  \\
25, 30, 35, 36, 47, 54, 55, 59, 67, 69, 87, 150, 195
\end{array} \right\}.
\]
This was calculated with the \textsc{Magma} function \texttt{IsIsomorphic}.
It is trivial that an automorphism of order $31$ of $D_{63,i}$ 
induces  an automorphism of order $31$ of $C'_2(D_{63,i})$.
Note that there are $38$ inequivalent extremal  doubly even self-dual codes of
length $64$ with  an automorphism of order $31$~\cite{Y}.

\begin{table}[thb]
\caption{Dimensions of $C'_2(D_{63,i})$}
\label{Tab:F2-1}
\centering
\medskip
{\small
\begin{tabular}{c|l}
\noalign{\hrule height1pt}
\multicolumn{1}{c|}{$\dim(C'_2(D_{63,i}))$ }& \multicolumn{1}{c}{$i$} \\
  \hline
$ 7$ &$805$\\
$12$& $806, 807, 808, 809, 810, 812, 813, 816, 820, 823$\\
$17$& $789, 790, 793, 800, 804$\\
$22$& $791, 795, 796, 797, 801, 802, 811, 814, 815, 817$, \\& $818, 819, 821, 822, 824, 826$\\
$32$& others\\ 
   \noalign{\hrule height1pt}
\end{tabular}
}
\end{table}

\begin{table}[thb]
\caption{Minimum weights of doubly even self-dual codes $C'_2(D_{63,i})$}
\label{Tab:F2-2}
\centering
\medskip
{\small
\begin{tabular}{c|l}
\noalign{\hrule height1pt}
\multicolumn{1}{c|}{$d$ }& \multicolumn{1}{c}{$i$} \\
  \hline
$4$ & $792$\\
$8$ &
$5, 14, 16, 17, 19, 21, 28, 30, 32, 34, 35, 36, 37, 41, 43, 49, 51, 53, 54, 59,$\\&
$ 60, 61, 64, 66, 67, 70, 73, 75, 78, 82, 85, 92, 93, 103, 107, 113, 114, 118, $\\&
$120, 123, 124, 135, 136, 138, 140, 141, 142, 143, 146, 153, 155, 156, 162, $\\&
$166, 171, 172, 173, 180, 182, 185, 186, 188, 189, 200, 205, 206, 208, 211, $\\&
$212, 213, 214, 229, 231, 232, 237, 242, 243, 244, 246, 248, 250, 251, 254, $\\&
$257, 259, 263, 264, 265, 266, 267, 270, 271, 272, 274, 275, 279, 280, 281, $\\&
$285, 287, 295, 297, 298, 300, 305, 309, 310, 313, 316, 319, 322, 323, 325, $\\&
$333, 337, 338, 341, 342, 343, 347, 354, 357, 358, 359, 369, 370, 381, 382, $\\&
$383, 390, 403, 407, 412, 413, 418, 420, 423, 425, 430, 432, 434, 441, 442, $\\&
$444, 446, 456, 458, 468, 471, 476, 477, 478, 483, 484, 485, 490, 497, 500, $\\&
$504, 505, 507, 509, 510, 515, 516, 519, 521, 527, 536, 537, 538, 546, 547, $\\&
$548, 553, 555, 558, 559, 565, 572, 575, 587, 588, 591, 594, 596, 598, 599, $\\&
$600, 601, 604, 606, 609, 610, 611, 612, 614, 617, 628, 634, 635, 636, 637, $\\&
$639, 642, 643, 645, 646, 647, 648, 651, 653, 654, 656, 657, 658, 659, 661, $\\&
$662, 663, 669, 672, 673, 674, 682, 686, 688, 689, 690, 697, 707, 711, 712, $\\&
$713, 715, 717, 718, 728, 730, 739, 743, 748, 752, 753, 756, 757, 758, 762, $\\&
$763, 764, 765, 766, 767, 768, 769, 770, 771, 772, 773, 775, 779, 780, 781, $\\&
$782, 783, 784, 785, 786, 787, 788, 798, 799, 825$
\\
$12$ & others \\   
         \noalign{\hrule height1pt}
\end{tabular}
}
\end{table}

\section{Summary}
\label{Sec:Summary}

We end this paper by listing a summary of the classification of 
Hadamard $2$-$(2p+1,p,(p-1)/2)$ designs with an automorphism of order $p$
having one fixed point
and the classification of Hadamard matrices of order $2p+2$ with an automorphism of order $p$
for $p \le 31$.
The number $N(\cD_{2p+1})$
of non-isomorphic Hadamard $2$-$(2p+1,p,(p-1)/2)$ designs 
with an automorphism of order $p$ having one fixed point
and the number $N(H_{2p+2})$ of
inequivalent Hadamard matrices of order $2p+2$ with an automorphism of order $p$
are listed in Table~\ref{Tab:Summary}, along with relevant references.

\begin{table}[thb]
\caption{Summary}
\label{Tab:Summary}
\centering
\medskip
{\small
\begin{tabular}{c|cc|c}
\noalign{\hrule height1pt}
$p$ & $N(\cD_{2p+1})$ & $N(H_{2p+2})$ & References\\
    \hline
 $3$ & $1$ & $1$ &\cite{Paley}\\
 $5$ & $1$ & $1$ &\cite{Paley} (see~\cite{HallB}) \\
 $7$ & $3$ & $3$  &\cite{Hall} \\   
 $11$ & $5$ & $4$     & \cite{H24} \\
$13$ & $7$ & $4$     &\cite{tonchev1}\\
$17$ & $21$ & $11$  &\cite{tonchev2}\\
$19$ & $33$ & $18$  &\cite{DHM}\\
$23$ & $109$ & $56$ &\cite{DHM}\\
$29$ & $531$ & $266$ &Section~\ref{Sec:60}\\
$31$ & $826$ & $414$ &Section~\ref{Sec:64}\\
\noalign{\hrule height1pt}
\end{tabular}
}
\end{table}

\bigskip
\noindent
\textbf{Acknowledgments.}
This work was supported by JSPS KAKENHI Grant Numbers 19H01802, 21K03350 and 23H01087.



     

\end{document}